\documentclass[11pt,draftcls,onecolumn,journal]{IEEEtran}
\usepackage{amsmath,amsfonts}
\usepackage{caption}      % for \captionof
\usepackage[justification=raggedright,singlelinecheck=false]{caption}
\usepackage{algorithmic}
\usepackage{algorithm}
\usepackage{array}
\usepackage[caption=false,font=normalsize,labelfont=sf,textfont=sf]{subfig}
\usepackage{textcomp}
\usepackage{stfloats}
\usepackage{url}
\usepackage{verbatim}
\usepackage{graphicx}
\usepackage{cite}
\usepackage{hyperref}
\usepackage{booktabs}
\usepackage{cite}
\usepackage{amsmath,amssymb,amsfonts}

\usepackage{graphicx}
\usepackage{textcomp}

\usepackage{bm}
\usepackage[most]{tcolorbox}
\tcbuselibrary{breakable}

\newtcolorbox{remarkbox}[2][]{%
  enhanced,                   % better title control
  breakable,
  colback=gray!5!white,       % box background
  colframe=gray!5!white,      % invisible frame
  boxrule=0pt,                 % no border line
  sharp corners,               % square corners
  colbacktitle=gray!10!white,  % title background
  fonttitle=\bfseries,         % bold title
  coltitle=blue!45!black,      % title text color
  title={#2},                  % title text
  before upper={%
  },
  after upper={%
  },
  #1                           % extra options
}

\usepackage{amsthm}
\theoremstyle{remark}

\def\BibTeX{{\rm B\kern-.05em{\sc i\kern-.025em b}\kern-.08em
    T\kern-.1667em\lower.7ex\hbox{E}\kern-.125emX}}

\usepackage[dvipsnames]{xcolor}
\usepackage{tikz}
\definecolor{accessblue}{cmyk}{1, 0.3, 0, 0.2}
\definecolor{greycolor}{cmyk}{0,0,0,.8}
\definecolor{revisionblue}{RGB}{0,0,0}%{0,70,140}
\usepackage{pgfplots}
\usepackage[nolist]{acronym}
\usepackage{booktabs} 		% cute tables
\usepackage{diagbox}		% diagonal cell division in table
\usepackage{upgreek}
\usepackage{bbm}
\usepackage{Files/winsnotation}
\usepackage{comment}
\usepackage{multirow}
\usepackage{mathtools}

\pgfplotsset{compat=1.18}

\begin{document}
%
% --------------------------------------------------------------------------
% Title
\title{A simple derivation of the Kalman filter %[Lecture Notes] 
}
%
% --------------------------------------------------------------------------
% Authors
\author{Marco~Chiani, %~\IEEEmembership{Fellow,~IEEE,}
Giovanni Petris, 
Moe Z. Win%,~\IEEEmembership{Fellow,~IEEE}, % <-this % stops a space
\thanks{MC is with the Department of Electrical, Electronic, and Information Engineering ``Guglielmo Marconi'', University of Bologna, Italy.}
\thanks{GP is with the Department of Mathematical Sciences, University of Arkansas, USA.}
\thanks{MZW is with the Massachusetts Institute of Technology, USA.}
\thanks{Manuscript received XXX 2025.}
}
\markboth{Chiani, Petris, Win, IEEE ,~Vol.~XX, No.~XX, ~2026}%
{Chiani \MakeLowercase{\textit{et al.}}: Kalman}


\date{}

\begin{acronym}
% usage: \ac{SW}, \acp{SW} for plurals \acf{SW} Use the full name of the acronym.
%\acs{SW}Use the acronym, even before the first corresponding \ac command
%\acl{acronym}Expand the acronym without using the acronym itself.
\small
\acro{AWGN}{additive white Gaussian noise}
\acro{BCH}{Bose–Chaudhuri–Hocquenghem}
\acro{BC}{bubble clustering}
\acro{BLUE}{best linear unbiased estimator}
\acro{CDF}{cumulative distribution function}
\acro{CRC}{cyclic redundancy code}
\acro{LDPC}{low-density parity-check}
\acro{LUT}{lookup table}
\acro{LSE}{least squares estimator}
\acro{ML}{maximum likelihood}
\acro{MWPM}{minimum weight perfect matching}
\acro{QECC}{quantum error correcting code}
\acro{PDF}{probability density function}
\acro{PMF}{probability mass function}
\acro{MPS}{matrix product state}
\acro{WEP}{weight enumerator polynomial}
\acro{WE}{weight enumerator}
\acro{BD}{bounded distance}
\acro{QLDPC}{quantum low density parity check}
\acro{CSS}{Calderbank, Shor, and Steane}
\acro{MST}{minimum spanning tree}
\acro{PruST}{pruned spanning tree}
\acro{RFire}{Rapid-Fire}
\acro{UF}{union-find}
\acro{LEMON}{library for efficient modeling and optimization in networks}
\acro{STM}{spanning tree matching}
\acro{i.i.d.}{independent identically distributed}
\acro{QEC}{quantum error correction}
\acro{BP}{belief propagation}

\end{acronym}

\maketitle
%
%
%
%
% ---------------------------------------------------
% Section
% ---------------------------------------------------
%
\section*{Scope}
\label{sec:abstract}

In this lecture note, we present a  concise and self-contained derivation of the discrete-time Kalman filter equations that requires only a basic understanding of least squares estimation. The treatment is designed to minimize mathematical overhead while preserving both rigor and generality.

\section*{Relevance}
The Kalman filter is one of the most important and widely used data fusion algorithms in information processing, with applications in statistical signal processing, radar tracking, inertial navigation, and economics \cite{Kal:60,KalFalArb:B69, Jaz:70,Kai:74, Har:90, Kay:93,KaiSayHas:B00, Sim:06, PetPetCam:09,GreAnd:15}.

Despite its importance, the derivation of the Kalman filter is often presented in a mathematically dense form that can be challenging for students and practitioners encountering it for the first time. As a result, many expositions omit the detailed steps, focusing instead on examples and applications.

This lecture note addresses that gap by providing a concise, self-contained, and pedagogically oriented derivation of the Kalman filter equations directly from the least squares estimation of the system state. The approach minimizes mathematical overhead, making the derivation accessible with only a basic understanding of least squares estimation, while preserving rigor and generality. The model and derivation do not assume Gaussian disturbances and can accommodate correlations between process and measurement errors. 

\subsection{Notation} Random variables are displayed in sans serif, upright fonts; their realizations in serif, italic fonts. Vectors are denoted by bold lowercase letters. Matrices are denoted by bold uppercase letters. For example, a random vector and its realization are denoted by $\RV{x}$ and $\V{x}$. 
The symbol $( )^\top$ is used for vector and matrix transposition. Symbol $\E{}$ represents the expectation operator. 
%
%
%
%
% ---------------------------------------------------
% Section
% ---------------------------------------------------
%
\section*{Prerequisites}
%\label{sec:notation}

We recall that, for a random vector $\RV{x}$ with $\E{\RV{x}}=\mathbf{0}$ (zero mean), the covariance matrix is  $\M{R}_{\RV{x}} =  \E{\RV{x} \RV{x}^\top}$. Given a deterministic matrix $\M{A}$, the random vector  $\RV{y}=\M{A} \RV{x}$ is zero mean with covariance  $\M{R}_{\RV{y}} =  \E{\RV{y} \RV{y}^\top} = \E{\M{A} \RV{x} \RV{x}^\top \M{A}^\top}=\M{A} \E{\RV{x} \RV{x}^\top} \M{A}^\top=\M{A} \M{R}_{\RV{x}} \M{A}^\top$.
%
%
%
%
% ---------------------------------------------------
% Section
% ---------------------------------------------------
%
%\section*{Tools}
%\label{sec:length}

We will need a few well-known tools from linear algebra. 
%
%\begin{itemize}
	%\item 
\begin{remarkbox}
{\textbf{Tool 1)}: Least Squares Estimator for the linear model: Gauss-Markov theorem} 

	Consider the general linear model \cite[p. 141]{Kay:93}
	\begin{equation}\label{eq:yan}
	\V{y}= \M{A} \V{x} + \V{n}
	\end{equation} 
	where $\V{y}$ is the measurements vector, the matrix $\M{A}$ is known and has full column rank, $\V{x}$ is a vector to be estimated,  and $\V{n}$ is a noise vector with zero mean and covariance matrix $\M{R}_{\RV{n}} =  \E{\RV{n} \RV{n}^\top}$. 
	
	The Gauss-Markov theorem states that, if $\M{R}_{\RV{n}} = \M{I}$ (white noise), the \ac{BLUE} is the ordinary \ac{LSE} of $\V{x}$, providing,\footnote{The result is easy to remember, as follows: if $\M{A}$ is a square matrix, we could find an unbiased estimate of $\V{x}$ as $\M{A}^{-1} \V{y}= \V{x} + \M{A}^{-1} \V{n}$. However, $\M{A}$ is in general a rectangular matrix, with more rows than columns, and we assume it is full rank. To build a square matrix we can first multiply by $\M{A}^{T}$ on the left, giving   $\M{A}^{T} \V{y}= \M{A}^{T} \M{A} \V{x} + \M{A}^{T} \V{n}$. Now, $\M{A}^{T} \M{A}$ is square and full rank, so we can use its inverse to get $\big( \M{A}^{T} \M{A} \big)^{-1} \M{A}^{T} \V{y} = \V{x} + \big( \M{A}^{T} \M{A} \big)^{-1} \M{A}^{T} \V{n} = \V{x}+\V{e}$. 
    The Gauss–Markov theorem states that, among all unbiased linear estimators, this estimator has the smallest error variance for each component, and therefore also minimizes the sum of all variances.
    } 	
	\begin{equation*}%\label{eq:xlsestimate}
	\widehat{\V{x}}= \big( \M{A}^\top \M{A} \big)^{-1}\M{A}^\top \V{y} = \V{x} + \V{e} \,.
	\end{equation*}
	In addition, the estimation error $\V{e}$ is zero mean with covariance 
	\begin{equation*}%\label{eq:coverrorafterestimate}
	\M{R}_{\RV{e}} = \E{\RV{e} \RV{e}^\top} = \big( \M{A}^\top \M{A} \big)^{-1} \,.
	\end{equation*}
In general, if the noise has covariance $\M{R}_{\RV{n}}$, the \ac{BLUE}  is the generalized \ac{LSE} of $\V{x}$ 	
	\begin{equation}\label{eq:xlsestimate}
	\widehat{\V{x}} = \big( \M{A}^\top \M{R}^{-1}_{\RV{n}} \M{A} \big)^{-1}\M{A}^\top \M{R}^{-1}_{\RV{n}} \V{y} = \V{x} + \V{e}
	\end{equation}
	where the estimation error $\V{e}$ is zero mean with covariance 
	\begin{equation}\label{eq:coverrorafterestimate}
	\M{R}_{\RV{e}} = \E{\RV{e} \RV{e}^\top} = \big( \M{A}^\top \M{R}^{-1}_{\RV{n}} \M{A} \big)^{-1} \,.
	\end{equation}
The latter expressions are derived from the former by first ``whitening'' the noise in \eqref{eq:yan} with  $\M{R}^{-1/2}_{\RV{n}} \V{y}= \M{R}^{-1/2}_{\RV{n}} \M{A} \V{x} + \M{R}^{-1/2}_{\RV{n}} \V{n} = \M{R}^{-1/2}_{\RV{n}} \M{A} \V{x} + \widetilde{\V{n}}$, where now   $\widetilde{\V{n}}$ has covariance $\E{\widetilde{\RV{n}}\widetilde{\RV{n}}^\top}=  \M{R}^{-1/2}_{\RV{n}}  \M{R}_{\RV{n}}  \M{R}^{-1/2}_{\RV{n}} = \M{I}$, and then applying the previous result.  	

	Note that the distribution of the noise is arbitrary. If it is Gaussian, then the \ac{LSE} is also the maximum likelihood estimator, and it is efficient as it achieves the Cramér--Rao bound \cite{Kay:93}.

\end{remarkbox}

\begin{remarkbox}{\textbf{Tool 2)} The Woodbury matrix inversion identity}
\begin{equation}\label{eq:woodbury}
\big( \M{A} + \M{U} \M{C} \M{V} \big)^{-1}
= \M{A}^{-1} - \M{A}^{-1} \M{U}
\big( \M{C}^{-1} + \M{V} \M{A}^{-1} \M{U} \big)^{-1}
\M{V} \M{A}^{-1}
\end{equation}	
\end{remarkbox}

\begin{remarkbox}{\textbf{Tool 3)} The matrix identity}
\begin{equation} \label{eq:wood2}
\M{C} \M{U}^\top \big( \M{A} + \M{U} \M{C} \M{U}^\top \big)^{-1}
= \big( \M{C}^{-1} \!+ \M{U}^\top \M{A}^{-1} \M{U} \big)^{-1}
\M{U}^\top \M{A}^{-1}
\end{equation} 
derived by writing $\big( \M{C}^{-1} + \M{U}^\top \M{A}^{-1} \M{U} \big) \,\M{C} \M{U}^\top = \M{U}^\top \M{A}^{-1} \big( \M{A} + \M{U} \M{C} \M{U}^\top \big)$. 
%\end{itemize}
\end{remarkbox}
%
%
% ---------------------------------------------------
% Section
% ---------------------------------------------------
%
\section*{Kalman Filter: Problem statement}
\label{sec:kalman}

The Kalman filter addresses the problem of estimating the internal state of a discrete-time linear dynamic system, given noisy observations and control inputs.  
%
%\subsection*{System Model}
%
The system is described by the state evolution and measurement equations:
\begin{align} \label{eq:modelx}
\V{x}' & = \M{F}\V{x}  + \M{B} \V{u} + \V{w}  \\
\V{z} &= \M{H} \V{x}'  + \V{v} \label{eq:modelz}
\end{align}
where $\V{x}' \in \mathbb{R}^s$ is the unknown current state, which depends linearly on the previous state $\V{x}$, and, possibly,  on a control input  $\V{u}$. We also know some measured data $\V{z}  \in \mathbb{R}^m $,  which depend on the current state $\V{x}'$. The measurement data $\V{z}$, typically representing a partial and noisy version of the state acquired through different sensors, can be used to improve the estimation of $\V{x}'$  (data fusion). 

We assume the following: 
\begin{itemize}
    \item The state evolution incorporates, in addition to the previous state and control vector $\V{u}$,  the process (or plant) noise $\V{w}$. Similarly, the  measurement equation incorporates the measurement (or observation) noise $ \V{v}$.   
    \item The process noise $\V{w}$ and observation noise $\V{v}$ are realizations of zero mean random vectors (the distribution is not relevant) with covariances $\M{Q}= \E{\RV{w}\RV{w}^\top}$ and $\M{R}= \E{\RV{v}\RV{v}^\top}$, respectively, and cross covariance $\M{C} =\E{\RV{w} \RV{v}^\top}$. 
    \item The matrices $\M{F}, \M{Q}, \M{H}, \M{R}$, as well the control $\M{B} \V{u}$, are perfectly known (the uncertainties are collectively modeled by  $\V{w}$ and $ \V{v}$). 
    \item We have an estimate $\widehat{\V{x}}=\V{x}+\V{e}$ of the previous state, where the estimation error $\V{e}$ is a realization of a zero mean random vector with known covariance $\M{P}=\E{\RV{e}\RV{e}^\top}$.  

\end{itemize}

\subsection*{Goal}

The problem for which we show the solution in this lecture notes article is how to find an estimate $\widehat{\V{x}}'=\V{x}'+\V{e}'$ of the current state, and the covariance $\M{P}'=\E{\RV{e}'(\RV{e}')^\top}$ of the corresponding estimation error  $\V{e}'$. 

\medskip

Together, equations~\eqref{eq:modelx} and~\eqref{eq:modelz} are often used as the basic building block to describe a system that evolves in discrete time. To make the dependence on time explicit, the equations are commonly written at time \( k \) as $
\V{x}_k = \M{F}_k \V{x}_{k-1} + \M{B}_k \V{u}_k + \V{w}_k$, $
\V{z}_k = \M{H}_k \V{x}_k + \V{v}_k$,  
and the relevant covariances as \( \M{Q}_k \), \( \M{R}_k \), \( \M{C}_k \), and \( \M{P}_k \). If $\E{\RV{w}_h^{\phantom{T}} \RV{w}_k^\top} = \M{0}$, $\E{\RV{v}_h^{\phantom{T}} \RV{v}_k^\top} = \M{0}$, and $\E{\RV{w}_h^{\phantom{T}} \RV{v}_k^\top} = \M{0}$, for $h\neq k$, then the Kalman filter can be applied iteratively, providing an estimate of the system state $\V{x}_k$ at the most recent time $k$ for which an observation $\V{z}_k$ is available. However, to keep the notation light, we will focus on an unspecified time point and drop the dependence on $k$, as we already did in~\eqref{eq:modelx} and~\eqref{eq:modelz}.

\section*{Solution: derivation of the Kalman equations}

We begin assuming that the process and observation noise are uncorrelated, i.e., $\E{\RV{w} \RV{v}^\top}=\mathbf{0}$.  
Then, we rewrite \eqref{eq:modelx} as 
\begin{align} \label{eq:systemmodel2}
\V{x}' & = \M{F} \big( \widehat{\V{x}}-\V{e} \big) + \M{B} \V{u} + \V{w}   = \M{F} \widehat{\V{x}} + \M{B} \V{u} + \V{w} -\M{F} \V{e}  
\end{align}

\noindent {\bf Key step:} from the system model \eqref{eq:modelx} and \eqref{eq:modelz}, moving all known quantities on the left we get
\begin{align} \label{eq:systemmodel3}
\M{F} \widehat{\V{x}} + \M{B} \V{u} &= \V{x}' +  \M{F} \V{e}-\V{w} \\
\V{z} &= \M{H} \V{x}'  + \V{v}
\end{align}
Then, using partitioned matrices, we write
\begin{equation}\label{eq:partitionedmodel}
\underbrace{\begin{bmatrix}
\M{F} \widehat{\V{x}} + \M{B} \V{u} \\
\V{z}
\end{bmatrix}
}_{\V{y}}= \underbrace{\begin{bmatrix}
 \M{I} \\
 \M{H}
\end{bmatrix}}_{\M{A}}
\V{x}'+
\underbrace{\begin{bmatrix}
 \M{F} \V{e}-\V{w} \\
 \V{v}
\end{bmatrix}}_{\V{n}}
\end{equation}
which is exactly in the form \eqref{eq:yan}, with 
\begin{equation} \label{eq:apart}
\V{y}= \begin{bmatrix}
\V{y}_1 \\
\V{y}_2
\end{bmatrix} \triangleq \begin{bmatrix}
\M{F} \widehat{\V{x}} + \M{B} \V{u} \\
\V{z}
\end{bmatrix} \qquad \M{A}\triangleq \begin{bmatrix}
 \M{I} \\
 \M{H}
\end{bmatrix} 
 \qquad  
\V{n}\triangleq \begin{bmatrix}
 \V{n}_1 \\
 \V{n}_2
\end{bmatrix}  
\end{equation}
where $\V{n}_2= \V{v}$, and $ \V{n}_1= \M{F} \V{e}-\V{w}$ is zero mean with covariance matrix 
\begin{equation}\label{eq:R1}
\widetilde{\M{Q}} \triangleq \E{ \RV{n}_1  \RV{n}_1^\top}=\M{F} \M{P} \M{F}^\top + \M{Q} \,.
\end{equation}
Thus, the covariance of $\RV{n}$ is 
\begin{equation} \label{eq:rpart}
\M{R}_{\RV{n}} = \E{\RV{n} \RV{n}^\top}= \E{ \begin{bmatrix}
 \RV{n}_1 \\
 \RV{n}_2
\end{bmatrix} \begin{bmatrix}
 \RV{n}^\top_1 &
 \RV{n}^\top_2
\end{bmatrix}}= \begin{bmatrix}
\widetilde{\M{Q}} & 0 \\
0 & \M{R} 
\end{bmatrix} \,.
\end{equation}

By applying the \ac{LSE} \eqref{eq:xlsestimate} and \eqref{eq:coverrorafterestimate} to \eqref{eq:partitionedmodel}, we have immediately the solution to our problem: 
\begin{align} \label{eq:solution1}\text{Kalman first form} \left\{ \begin{array}{lcl}
\widehat{\V{x}}' & = & \big(\M{A}^\top \M{R}^{-1}_{\RV{n}} \M{A} \big)^{-1}\M{A}^\top \M{R}^{-1}_{\RV{n}} \V{y} \\
 \M{P}' & = & \big(\M{A}^\top \M{R}^{-1}_{\RV{n}} \M{A} \big)^{-1}
\end{array}\right. \,. 
\end{align} 

The equations \eqref{eq:solution1} give the \ac{BLUE} of the state $\V{x}'$, taking into account also the observed data, and are equivalent to the standard form of the Kalman filter equations, which is developed next. Note that $\M{R}_{\RV{n}}$ has dimensions $(s+m)\times (s+m)$.

Due to the diagonal form of $\M{R}_{\RV{n}}$ in \eqref{eq:rpart}, we may use some matrix identities to work with smaller matrices, obtaining the standard form of the Kalman filter. 
First, with \eqref{eq:apart} and \eqref{eq:rpart} we have 

$$\M{A}^\top \M{R}^{-1}_{\RV{n}} \M{A} = \begin{bmatrix}
\M{I} & \M{H}^\top
\end{bmatrix}  \begin{bmatrix}
\widetilde{\M{Q}}^{-1}  & 0 \\
0 & \M{R}^{-1}
\end{bmatrix} \begin{bmatrix}
 \M{I} \\
 \M{H}
\end{bmatrix} = \widetilde{\M{Q}}^{-1} +  \M{H}^\top \M{R}^{-1} \M{H} $$  
so that we can write \eqref{eq:solution1} in a second form 
\begin{align} \label{eq:solution2}\text{Kalman 2nd form} \left\{ \begin{array}{lcl}
\widehat{\V{x}}' & = & \big( \widetilde{\M{Q}}^{-1} + \M{H}^\top \M{R}^{-1} \M{H} \big)^{-1} \big( \widetilde{\M{Q}}^{-1} \big(\M{F} \widehat{\V{x}} + \M{B} \V{u} \big) + \M{H}^\top \M{R}^{-1} \V{z} \big)\\
\M{P}'  & = & \big(\widetilde{\M{Q}}^{-1} + \M{H}^\top \M{R}^{-1} \M{H} \big)^{-1} 
\end{array}\right. \,. 
\end{align}
Finally, applying $\textbf{Tool 2}$ and $\textbf{Tool 3}$ to respectively the first and second term in \eqref{eq:solution2}, we have the standard form of the Kalman filter equations
\begin{align} \label{eq:solution3}\text{} \left\{ \begin{array}{lcl}
\widehat{\V{x}}' & = & \big( \M{I}- \widetilde{\M{Q}}  \M{H}^\top \big(\M{R} +  \M{H} \widetilde{\M{Q}}  \M{H}^\top \big)^{-1}  \M{H} \big) \big(\M{F} \widehat{\V{x}} + \M{B} \V{u} \big) 
+ \widetilde{\M{Q}} \M{H}^\top \big(\M{R} +  \M{H} \widetilde{\M{Q}}  \M{H}^\top \big)^{-1} \V{z} \\ %\RV{z} \\
\M{P}'  & = & \big( \M{I}- \widetilde{\M{Q}}  \M{H}^\top \big( \M{R} + \M{H} \widetilde{\M{Q}} \M{H}^\top \big)^{-1}  \M{H} \big) \,\widetilde{\M{Q}}  
\end{array}\right. 
\end{align}
that are equivalent to \eqref{eq:solution1} and \eqref{eq:solution2}, but with the advantage of requiring just one $m \times m$ matrix inversion. 
By denoting the Kalman gain 
\begin{equation} \label{eq:kgain} 
\M{K} = \widetilde{\M{Q}} \M{H}^\top \big( \M{R} + \M{H} \widetilde{\M{Q}}  \M{H}^\top \big)^{-1} 
\end{equation}
we can rewrite expression \eqref{eq:solution3} as
\begin{align} \label{eq:solution4}\text{Kalman 3rd form} \left\{ \begin{array}{lcl}
\widehat{\V{x}}' & = & \big(\M{I}-  \M{K}  \M{H} \big) \big( \M{F} \widehat{\V{x}} + \M{B} \V{u} \big) + \M{K} \V{z}\\
\M{P}' & = & \big(\M{I}- \M{K}  \M{H} \big) \,\widetilde{\M{Q}}  
\end{array}\right. 
\end{align}
where we recall that $\widetilde{\M{Q}}=\M{F} \M{P} \M{F}^\top + \M{Q}$. 

The equations \eqref{eq:solution4}, together with the gain \eqref{eq:kgain}, express the Kalman filter in the standard form.

\begin{remarkbox}{Generalization: correlated  process and measurement noise}
In the general case where the process noise and the measurement noise have correlation $\M{C} =\E{\RV{w} \RV{v}^\top}$, recalling the definition $ \RV{n}_1= \M{F} \RV{e}-\RV{w}$, the covariance matrix in \eqref{eq:rpart} becomes 
\begin{equation} \label{eq:rpartcorrelated}
\M{R}_{\RV{n}}= \E{\RV{n} \RV{n}^\top}= \begin{bmatrix}
\widetilde{\M{Q}} & -\M{C} \\
-\M{C}^\top & \M{R} 
\end{bmatrix} 
\end{equation}
which is no longer diagonal. We can still use \eqref{eq:solution1}, but if we want to spare inversions we will need to revise the simplifications.  An easy way is to make a linear transformation on the $\V{y}$ to decorrelate $\RV{n}_1$ and $\RV{n}_2$, to then use the previous results.
Let us look at what happens for the transformation  
$\RV{n}_2 \rightarrow \widetilde{\RV{n}}_2=\RV{n}_2 + \M{D} \RV{n}_1$, where $\M{D}$ is a suitable matrix. 
In this case $\E{\RV{n}_1 \widetilde{\RV{n}}_2^\top}=-\M{C}+ \widetilde{\M{Q}} \M{D}^\top$. 
Thus, by chosing $\M{D}=\M{C}^\top\widetilde{\M{Q}}^{-1}$, we can decorrelate the noise parts $\RV{n}_1$ and $\widetilde{\RV{n}}_2$. 

In other words, we can decorrelate the noise terms multiplying  \eqref{eq:partitionedmodel} by the matrix
\begin{equation} \label{eq:Gdecorrelate}
\mathbf{G}=
\begin{bmatrix}
\M{I} & \mathbf{0} \\
\M{C}^\top \widetilde{\M{Q}}^{-1} & \M{I} 
\end{bmatrix} 
\end{equation}
which produces a noise vector with covariance\footnote{We can see $\mathbf{G}^\top$ as a congruence transformation which diagonalizes $\mathbf{G} \M{R}_{\RV{n}} \mathbf{G}^\top$.}  
\begin{equation} \label{eq:Rdecorrelated}
\mathbf{G} \M{R}_{\RV{n}} \mathbf{G}^\top=
\begin{bmatrix}
\widetilde{\M{Q}} & \mathbf{0} \\
\mathbf{0} & \M{R}-\M{C}^\top \widetilde{\M{Q}}^{-1} \M{C} 
\end{bmatrix} \,.
\end{equation}

Therefore, 
the decorrelated problem is like \eqref{eq:partitionedmodel} with the substitutions 
$\V{y}_2 \rightarrow %\widetilde{\V{y}}_2=
\V{y}_2+\M{C}^\top\widetilde{\M{Q}}^{-1}\V{y}_1$,   
$\M{H} \rightarrow %\widetilde{\M{H}}=
\M{H}+\M{C}^\top\widetilde{\M{Q}}^{-1}$, and correlation matrix \eqref{eq:Rdecorrelated}. 
With these changes in \eqref{eq:kgain} and \eqref{eq:solution4}, 
the Kalman gain in the general case becomes
\begin{equation} \label{eq:kgaincorr} 
  \M{K} = \big(\widetilde{\M{Q}} \M{H}^\top + \M{C} \big) \big( \M{R} + \M{H} \M{C} + \M{C}^\top \M{H}^\top + \M{H} \widetilde{\M{Q}} \M{H}^\top \big)^{-1} 
\end{equation}
and the Kalman filter equations become 
\begin{align} \label{eq:solution5}\text{Kalman 4th form} \left\{ \begin{array}{lcl}
\widehat{\V{x}}' & = & \big(\M{I}-  \M{K}  \M{H} \big) \big(\M{F} \widehat{\V{x}} + \M{B} \V{u} \big) +  \M{K} \V{z}\\
\M{P}'  & = & \widetilde{\M{Q}}- \M{K} \big(\M{H} \widetilde{\M{Q}} +  \M{C}^\top \big)  
\end{array}\right. 
\end{align}
where $\widetilde{\M{Q}}=\M{F} \M{P} \M{F}^\top + \M{Q}$.

Expressions \eqref{eq:kgaincorr} and \eqref{eq:solution5} are the Kalman filter equations, valid also when there is correlation between the process and the measurements errors. We note that, apart from the change in the Kalman gain, the state estimator has exactly the same form as for the uncorrelated case \eqref{eq:solution4}. 
\end{remarkbox}

\section*{Interpretation}

 The Kalman filter equations lend themselves to a variety of interpretations and applications, many of which are discussed in standard textbooks on Kalman filtering. Here, we highlight only a few points and refer the reader to the bibliography for a complete treatment~\cite{Jaz:70,Sim:06,Kay:93,Har:90,PetPetCam:09,GreAnd:15}.

\begin{itemize}
\item The derivation clearly shows that the Kalman filter is the least squares estimator of the state and, as such, is also the Maximum Likelihood estimator in the special case of normally distributed noise. Within a Bayesian framework, assuming that also $\RV x$ (or equivalently $\RV e$), is normally distributed, then $\RV x'$ is also normally distributed. In this case $\widehat{\V{x}}'$ and $\M{P}'$ can be interpreted as the conditional mean and conditional covariance of $\RV{x}'$, given the observation $\RV{z} = \V{z}$. Therefore, $\widehat{\V{x}}'$ is the MSE-optimal point estimate of $\V{x}'$.

\item Note that a prediction based solely on the state equation~\eqref{eq:modelx} is $\widehat{\V{x}}_\text{no-data}' \triangleq \M{F} \widehat{\V{x}} + \M{B} \V{u}$. Then, the estimate \eqref{eq:solution4} can be interpreted as 
\begin{align}\label{eq:KFxnodata}
\widehat{\V{x}}' = \widehat{\V{x}}_\text{no-data}' + \M{K} \big( \V{z} - \M{H} \widehat{\V{x}}_\text{no-data}' \big)    
\end{align}  
that is, the no-data prediction is corrected using the observed measurement through the Kalman gain~\( \M{K} \), applied to the residual \( \V{z} - \M{H} \widehat{\V{x}}_\text{no-data}' \), commonly referred to as the \emph{innovation}. Thus, the estimation process can be understood as a no-data prediction step followed by a data-driven correction step.  

\item The matrix  
\[
\widetilde{\M{Q}} = \M{F} \M{P} \M{F}^\top + \M{Q}
\]  
represents the error covariance obtained when estimating \( \V{x}' \) using only the state equation~\eqref{eq:modelx}. For this reason, it can be interpreted as \( \M{P}'_\text{no-data} \).  

%{\color{blue}
\item In the case of uncorrelated prediction and observation errors, the
Kalman estimate in~\eqref{eq:solution2} can also be characterized
variationally as
\begin{equation}
\label{eq:variationalKalman}
\widehat{\V{x}}'
=
\underset{\V{\xi}}{\operatorname{arg\,min}}
\left\{
(\V{\xi}-\widehat{\V{x}}'_{\mathrm{no-data}})^\top
\widetilde{\M{Q}}^{-1}
(\V{\xi}-\widehat{\V{x}}'_{\mathrm{no-data}})
+
(\V{z}-\M{H}\V{\xi})^\top
\M{R}^{-1}
(\V{z}-\M{H}\V{\xi})
\right\}.
\end{equation}
The two terms measure, respectively, the discrepancy from the model-based
prediction and the observation residual, weighting them according to their
uncertainties. Thus, information affected by greater uncertainty has less
influence on the resulting estimate. It is important to note, however, that
the weighting matrices in~\eqref{eq:variationalKalman} are not chosen
arbitrarily: their identification with the inverses of the prediction- and
observation-error covariances follows here from the Gauss--Markov formulation
of the stacked linear model. Setting the gradient of the functional to zero
gives exactly~\eqref{eq:solution2}. Hence, under the second-order assumptions
used here and without requiring Gaussian distributions, the variational form
provides an alternative characterization of the estimate already derived.
%}

\item We observe that \eqref{eq:solution1} requires the inversion of matrices of dimension $(s+m)\times(s+m)$, \eqref{eq:solution2} requires the inversion of matrices of dimensions $m\times m$ and $s\times s$, and \eqref{eq:solution4} requires the inversion of a matrix of dimension $m\times m$. Thus, if there is only a scalar measurement ($m=1$), no matrix inversion is needed. A method to avoid matrix inversion for uncorrelated measurements (sequential Kalman filtering) is discussed in the box below.

{\color{revisionblue}
\item We remark that equations~\eqref{eq:modelx}--\eqref{eq:modelz} are typically applied to follow the system evolution over discrete time steps $k = 0, 1, . . .$. 
A key consequence of the Kalman equations is that, at time $k$, to compute the new estimate $\widehat{\V{x}}_k$ it is sufficient to know the estimate and covariance at the previous step $\big( \widehat{\V{x}}_{k-1},\M{P}_{k-1} \big)$ and the current measurement $\V{z}_k$, given the model matrices at time $k$. More details on this issue are given in a following section.  

\item In practical implementations, the initial state covariance $\M{P}_0$ is often chosen as a diagonal matrix with large entries, reflecting high initial uncertainty.  
A common choice is $\M{P}_0 = \sigma_0^2 \M{I}$ with a large $\sigma_0^2$ when no reliable prior information is available.  
The process noise covariance $\M{Q}$ and the measurement noise covariance $\M{R}$ are usually tuned empirically.  
The matrix $\M{Q}$ controls the filter’s responsiveness to changes in the system dynamics, while $\M{R}$ governs the sensitivity to measurement noise.  
Underestimating $\M{Q}$ or overestimating $\M{R}$ may result in slow convergence or filter divergence, whereas an excessively large $\M{Q}$ can cause overly-reactive and noisy estimates.  
In practice, $\M{Q}$ and $\M{R}$ are often adjusted iteratively based on innovation (residual) statistics, performance criteria, or prior system knowledge.
}
\end{itemize}

\begin{remarkbox}{Sequential Kalman filtering}
From \eqref{eq:kgain} or \eqref{eq:kgaincorr}, computing the Kalman gain requires inverting an \( m \times m \) matrix. If \( m = 1 \) (a scalar measurement), this inversion reduces to a scalar division. More generally, if the process and measurement noise are uncorrelated and the \( m \) measurements are mutually uncorrelated (e.g., from \( m \) independent sensors), so that \( \M{R} \) is diagonal, matrix inversions can be avoided entirely by updating the estimates sequentially between time steps. For the first scalar sensor, one applies the Kalman equations with~\eqref{eq:modelx}; for the remaining \( m-1 \) sensors, one sets \( \M{F}=\M{I} \) and \( \M{Q}=\M{B}=\mathbf{0} \). This approach is known as \emph{sequential Kalman filtering}~\cite{Sim:06} (see Algorithm~1). The same principle applies if \( \M{R} \) is block diagonal, in which case the update is performed sequentially for each block of correlated sensors.

\bigskip 

\noindent\captionof{algorithm}{Sequential Kalman Estimator. Assumes $\M{C}=\M{0}$, measurements $\V{z} = ( z_1, z_2, \ldots,z_m )$, $z_i= \V{h}_i^\top \V{x}' + v_i$,  $\M{R} = \mathrm{diag} \big( \sigma^2_1, \sigma^2_2, \dots, \sigma^2_m \big)$}
\begin{algorithmic}[]
         
\STATE Initialize: $\V{x}_0 \gets \M{F}\widehat{\V{x}}+\M{B}\V{u}$, \quad $\widetilde{\M{Q}} \gets \M{F} \M{P} \M{F}^\top + \M{Q}$
\FOR{$i = 1$ \TO $m$}
%  \STATE Extract row: $\V{h}_i^\top \gets$ $i$-th row of $\M{H}$
  \STATE Kalman gain: $\M{K} \gets \widetilde{\M{Q}} \V{h}_i / \big( \sigma^2_i + \V{h}_i^\top \widetilde{\M{Q}} \V{h}_i \big)$
  \STATE State update: $\V{x}_i \gets \big( \M{I} - \M{K} \V{h}_i^\top \big) \V{x}_{i-1} + \M{K} z_i
$
  \STATE Covariance update: $\widetilde{\M{Q}} \gets \big(\M{I}-  \M{K}  \V{h}_i^\top \big) \,\widetilde{\M{Q}}$
\ENDFOR
\RETURN $\widehat{\V{x}}' \gets \V{x}_m$, \quad $\M{P}' \gets \widetilde{\M{Q}}$
\end{algorithmic}
%\end{algorithm}
\end{remarkbox}

{\color{revisionblue}
\section*{Worked numerical example (Kalman filtering)}

To illustrate the Kalman equations derived above, we consider a simple numerical example of a two-dimensional linear system describing the motion of a robot along a straight line. The initial speed is about $4\,\mathrm{m/s}\approx 14.4\,\mathrm{km/h}$, and position measurements are available every $\Delta t = 0.5\,\mathrm{s}$. For simplicity, physical units are omitted in the following.

The model follows~\eqref{eq:modelx}--\eqref{eq:modelz}, where the state vector is
\[
\V{x}=\begin{bmatrix} d \\ \dot d \end{bmatrix}
\]
with $d$ denoting position and $\dot d$ denoting velocity. The measurement is scalar ($m=1$) and provides information only on the position component (not on the velocity).

We use a constant-velocity nominal model over one sampling interval $\Delta t$:
\[
d\,' = d + \dot d\,\Delta t,\qquad \dot d\,'=\dot d,
\qquad\Rightarrow\qquad
\M{F}=
\begin{bmatrix}
1 & \Delta t\\
0 & 1
\end{bmatrix}
=
\begin{bmatrix}
1 & 0.5\\
0 & 1
\end{bmatrix}.
\]
For the state evolution we use~\eqref{eq:modelx} with $\M{B}\V{u}=\V{0}$. We assume a process noise covariance $\M{Q}=0.01\,\M{I}$. The measurement equation is~\eqref{eq:modelz} with $\M{H}=\big[1 \,\, 0\big]$ and measurement-noise covariance $\M{R}=0.04$.

Assume that the initial estimate and covariance at $k=0$ are
\[
\widehat{\V{x}}_0=
\begin{bmatrix}
\widehat d_0\\ \widehat{\dot d}_0
\end{bmatrix}
=
\begin{bmatrix}
0\\ 4
\end{bmatrix}, \qquad
\M{P}_0=
\begin{bmatrix}
1 & 0\\
0 & 1
\end{bmatrix}
\]
and that the first two measurements are $z_1=1.8$ and $z_2=4.3$.

\paragraph{Step 1} ($k=1$, measurement $z_1=1.8$)

\emph{No-data estimate (prediction)}
\[
\M{F}\widehat{\V{x}}_0=
\begin{bmatrix}2\\4\end{bmatrix},\qquad
%\M{P}_{1,\text{no-data}}=
\widetilde{\M{Q}}_1=\M{F}\M{P}_0\M{F}^\top+\M{Q}
=
\begin{bmatrix}1.26&0.50\\0.50&1.01\end{bmatrix}.
\]

\emph{Kalman gain and innovation}
\[
\M{K}_1 = \widetilde{\M{Q}}_1\M{H}^\top \big( \M{R} + \M{H} \widetilde{\M{Q}}_1 \M{H}^\top \big)^{-1}
=
\begin{bmatrix}0.969\\0.385\end{bmatrix},\quad
\nu_1=z_1-\M{H}\M{F}\widehat{\V{x}}_0=1.8-2.0=-0.2.
\]

\emph{Data-fusion step (measurement update)}
\[
\widehat{\V{x}}_1 = \M{F}\widehat{\V{x}}_0+\M{K}_1\nu_1=
\begin{bmatrix}1.806\\3.923\end{bmatrix},\qquad
\M{P}_1 = \big( \M{I} - \M{K}_1\M{H} \big) \,\widetilde{\M{Q}}_1
=
\begin{bmatrix}0.0388&0.0154\\0.0154&0.8177\end{bmatrix}.
\]

\paragraph{Step 2} ($k=2$, measurement $z_2=4.3$)

\emph{No-data estimate (prediction)}
\[
%\widehat{\V{x}}_{2,\text{no-data}}=
\M{F}\widehat{\V{x}}_1 =
\begin{bmatrix}3.768\\3.923\end{bmatrix},\qquad
%\M{P}_{2,\text{no-data}}=
\widetilde{\M{Q}}_2=\M{F}\M{P}_1\M{F}^\top+\M{Q}
=
\begin{bmatrix}0.2686&0.4242\\0.4242&0.8277\end{bmatrix}.
\]

\emph{Kalman gain and innovation}
\[
\M{K}_2=\widetilde{\M{Q}}_2\M{H}^\top \big( \M{R} + \M{H} \widetilde{\M{Q}}_2 \M{H}^\top \big)^{-1}
=
\begin{bmatrix}0.870\\1.375\end{bmatrix},\quad
\nu_2=z_2-\M{H}\M{F}\widehat{\V{x}}_1=4.3-3.768=0.532.
\]

\emph{Data-fusion step (measurement update)}
\[
\widehat{\V{x}}_2=\M{F}\widehat{\V{x}}_1+\M{K}_2\nu_2=
\begin{bmatrix}4.231\\4.655\end{bmatrix},\qquad
\M{P}_2 = \big( \M{I} - \M{K}_2 \M{H} \big)\,\widetilde{\M{Q}}_2
=
\begin{bmatrix}0.0348&0.0550\\0.0550&0.2445\end{bmatrix}.
\]

\medskip
\noindent\emph{Summary:}\;
$$\widehat{\V{x}}_0=\begin{bmatrix}0\\4\end{bmatrix}
\ \xrightarrow{\,z_1=1.8\,}\
\widehat{\V{x}}_1=\begin{bmatrix}1.806\\3.923\end{bmatrix}
\ \xrightarrow{\,z_2=4.3\,}\
\widehat{\V{x}}_2=\begin{bmatrix}4.231\\4.655\end{bmatrix}.$$

\medskip
\noindent\emph{Comments.}
 Here, only the position is directly observed; the velocity is inferred through the dynamics encoded in $\M{F}$. The example illustrates how the Kalman filter combines the model-based prediction $\widehat{\V{x}}_{k,\text{no-data}}$ and the noisy measurement $z_k$ to improve both the estimated position and the inferred velocity.  
Because $m=1$, the innovation covariance
$\M{R}+\M{H}\widetilde{\M{Q}}_k\M{H}^\top$
is a scalar, so the gain computation requires only scalar divisions (no matrix inversion).  
Finally, the decrease of $\M{P}_k$ from $\M{P}_0$ to $\M{P}_2$ reflects the increased confidence after incorporating measurements.

}

\newpage 

{\color{revisionblue}
%-------------------------------------------------------
\section*{Sequence of observations and Kalman: filtering and smoothing}
%-------------------------------------------------------

In the previous section we derived, for a \emph{single} prediction--update step,
the least--squares/BLUE estimate of the current state $\V x'$ given $(\widehat{\V x},\M P)$ and a measurement $\V z$,
obtaining the equivalent forms \eqref{eq:solution1}--\eqref{eq:solution4} and the gain \eqref{eq:kgain}.
We now show how the \emph{iterative} application of those equations yields the optimal least-squares/BLUE estimate
of $\V x_k$ given $\V z_{1:k}\triangleq(\V z_1, \V z_2, \dots, \V z_k)$ and how to compute the optimal \emph{smoothed} estimate
of $\V x_k$ by incorporating future measurements, i.e.,\ given $\V z_{1:n}$ with $n>k$.

\subsection*{Model and estimation tasks}
Consider the discrete-time linear model, for $k = 1, 2, \dots,n$,
\begin{align*}
\V x_k & = \M F \V x_{k-1}+\V w_k\\
\V z_k & = \M H \V x_k+\V v_k
\end{align*}
with zero-mean noises and covariances
\[
\E{\RV w_k\RV w_k^\top}=\M Q,\qquad \E{\RV v_k\RV v_k^\top}=\M R,
\qquad \E{\RV w_k\RV v_\ell^\top}=\M 0
\]
and all noise terms mutually independent across time.
Assume an unbiased prior $(\widehat{\V x}_0,\M P_0)$.

% %\paragraph{Filtering vs.\ smoothing.}
The dynamical model leads naturally to different estimation tasks:
\begin{enumerate}
  \item Estimate the final state $\V x_n$ using all measurements available at time $n$, i.e., $\V z_{1:n}$. We will indicate this estimate as $\widehat{\V x}_{n|n}$. 
  \item Estimate an intermediate state $\V x_k$ (or the whole trajectory $\V x_{1:n}$) using all measurements up to time $n$, i.e., $\V z_{1:n}$. We will indicate this estimate as $\widehat{\V x}_{k | n}$.
\end{enumerate}
We will show that Problem~1 can be solved recursively by the classical \emph{Kalman filtering}.
Problem~2 can be solved in \emph{batch mode} (wait to have all measurements and then solve for the whole state sequence), or by first solving recursively Problem 1 and then performing a backward recursion, obtaining the smoothed estimate of $\V x_k$ and the smoothing variance.  
% in what is called  \emph{smoothing}. 
Note that the estimate of $\V x_k$ that exploits all measurements $\V z_{1:n}$ is, in general, different from the filtered estimate of $\V x_k$ that uses only $\V z_{1:k}$.

\medskip
\noindent

To address the problem we proceed again with least-square estimation as in previous sections. About the state evolution, writing $\V x_0=\widehat{\V x}_0-\V e_0$ and moving known quantities to the left gives
\begin{align*}
-\M F\widehat{\V x}_0 &= -\V x_1+\widetilde{\V w}_1,\\
\M 0 &= \M F\V x_{k-1}-\V x_k+\V w_k,\qquad k=2,\dots,n,
\end{align*}
where $\widetilde{\V w}_1\triangleq \V w_1-\M F\V e_0$ has covariance $\E{\widetilde{\RV w}_1\widetilde{\RV w}_1^\top}=\M Q+\M F\M P_0\M F^\top\;\triangleq\;\M{\widetilde Q}_1$.

%-------------------------------------------------------
\subsection*{Batch (whole-sequence) least squares}
%-------------------------------------------------------

Stacking the equations for $k=1,\dots,n$ yields, similarly to \eqref{eq:partitionedmodel}, the linear model $\V{y}= \M{A} \V{x} + \V{n}$ as

\begin{equation}\label{eq:batchmodel12_improved}
\underbrace{\begin{bmatrix}
-\M F\widehat{\V x}_0\\
\V z_1 \\ 
\M 0\\
\V z_2\\ 
\vdots\\
\M 0\\
\V z_n \\
\end{bmatrix}}_{\V y}
=
\underbrace{\begin{bmatrix}
-\M I & \M 0 & \M 0 & \cdots & \M 0\\
\M H  & \M 0 & \M 0 & \cdots & \M 0\\ 
\M F  & -\M I& \M 0 & \cdots & \M 0\\
\M 0  & \M H & \M 0 & \cdots & \M 0\\ 
\vdots     & \vdots    & \ddots    & \ddots & \vdots\\
\M 0  & \M 0 & \cdots    & \M F & -\M I\\
\M 0  & \M 0 & \cdots    & \M 0 & \M H
\end{bmatrix}}_{\M A}
\underbrace{\begin{bmatrix}
\V x_1\\
\V x_2\\
\vdots\\
\V x_n
\end{bmatrix}}_{\V x}
+
\underbrace{\begin{bmatrix}
\widetilde{\V w}_1\\
\V v_1\\ 
\V w_2\\
\V v_2\\ 
\vdots\\
\V w_n\\
\V v_n
\end{bmatrix}}_{\V n}.
\end{equation}

The least-squares estimate $\widehat{\V x}$ of the whole vector $\V x$ satisfies the normal equation
\[
\M A^\top \M R_{\RV n}^{-1}\M A\,\widehat{\V x}
=
\M A^\top \M R_{\RV n}^{-1}\V y
\]
which, considering that the noise covariance $\M R_{\RV n} =\E{\RV n\RV n^\top}$ is block diagonal, and the structure of $\M A$ in \eqref{eq:batchmodel12_improved}, takes the block tridiagonal form $\M M \widehat{\V x} = \V b$ as follows 
\begin{equation}\label{eq:block_tridiag_improved}
\underbrace{\begin{bmatrix}
\M D_1 & \M E^\top & \M 0 & \cdots & \M 0\\
\M E & \M D_2 & \M E^\top & \cdots & \M 0\\
\M 0 & \M E & \M D_3 & \ddots & \vdots\\
\vdots & \ddots & \ddots & \ddots & \M E^\top\\
\M 0 & \cdots & \M 0 & \M E & \M D_n
\end{bmatrix}}_{\M M}
\underbrace{\begin{bmatrix}
\widehat{\V x}_{1|n}\\
\widehat{\V x}_{2|n}\\
\vdots\\
\widehat{\V x}_{n|n}
\end{bmatrix}}_{\widehat{\V x}}
=
\underbrace{\begin{bmatrix}
\V b_1\\
\V b_2\\
\vdots\\
\V b_n
\end{bmatrix}}_{\V b}
\end{equation}
where %Thus one obtains
\begin{equation}
\label{eq:Eb1bkD}
\begin{alignedat}{2}
\M E   %&= \M A_0^\top \M C^{-1}\M B_0
      &= -\,\M Q^{-1}\M F,\\
\V b_1 %&= \M A_0^\top \M C_1^{-1}\V y_1+\M B_0^\top\M C^{-1}\V y_2
      &= \M{\widetilde Q}_1^{-1}\M F\widehat{\V x}_0+\M H^\top\M R^{-1}\V z_1,\\
\V b_k %&= \M A_0^\top \M C^{-1}\V y_k
      &= \M H^\top\M R^{-1}\V z_k,\qquad k=2,\dots,n \\
    \M D_1 &= \M{\widetilde Q}_1^{-1}+\M H^\top \M R^{-1}\M H +\M F^\top \M Q^{-1}\M F \\
\M D_k & = \M Q^{-1}+\M H^\top \M R^{-1}\M H +\M F^\top \M Q^{-1}\M F,\qquad k=2,\dots,n-1 \\
\M D_n &= \M Q^{-1}+\M H^\top \M R^{-1}\M H \,.
\end{alignedat}
\end{equation}

Solving the system \eqref{eq:block_tridiag_improved}, for example by matrix inversion $\widehat{\V x}= \M M^{-1} \V b$, will give in a batch the whole vector $\widehat{\V x}$, i.e., the least-square estimation of all states $\widehat{\V x}_{1:n}$ considering all observations. However, we can avoid inverting the big matrix $\M M$ by solving the system via Gaussian elimination,  
%-------------------------------------------------------
%\subsection*{Recursive solution (block elimination) and its connection to Kalman filtering}
%-------------------------------------------------------
%
%A standard way to solve \eqref{eq:block_tridiag_improved} is by  Gaussian elimination, 
which transforms the matrix in an upper triangular form. It consists of a forward recursion, which at the end allows to determine $\widehat{\V x}_{n|n}$, eventually followed by a backward recursion to find all other $\widehat{\V x}_{k|n}$. For the linear system of equations \eqref{eq:block_tridiag_improved}, Gaussian elimination is achieved as follows.  

\noindent{\bf Forward recursion:} 

Initialize
\[
\bar{\M D}_1=\M D_1,\qquad \bar{\V b}_1=\V b_1
\]

and for $k=2,\dots,n$ compute
\begin{equation}
\label{eq:Dkbk}
\begin{alignedat}{2}
\bar{\M D}_k &= \M D_k-\M E\,\bar{\M D}_{k-1}^{-1}\M E^\top,\\
\bar{\V b}_k &= \V b_k-\M E\,\bar{\M D}_{k-1}^{-1}\bar{\V b}_{k-1}.
\end{alignedat}
\end{equation}

After the forward sweep, the system is transformed into upper triangular form, so 
\begin{align}\label{eq:xnnGeneral}
\widehat{\V x}_{n|n}=\bar{\M D}_n^{-1}\,\bar{\V b}_n\,.
\end{align}
\bigskip 

\noindent
For Problem 1,  we can stop here. 
For Problem 2, once we have $\widehat{\V x}_{n|n}$ we could also perform a backward substitution to recover all other $\widehat{\V x}_{k|n}$. 
In fact, after the forward recursion the system has been transformed into an upper block-bidiagonal one,
with blocks $\bar{\M D}_{k}$ on the main diagonal and $\M E^\top$ on the first upper diagonal. Therefore, the remaining unknowns
are recovered by the backward (back-substitution) recursion.  

\noindent
{\bf Backward recursion:} 
For $k=n-1,n-2,\ldots,1$ compute
\begin{equation}\label{eq:backwardgeneral}
\widehat{\V x}_{k|n} = \bar{\M D}_k^{-1} \!\big( \bar{\V b}_k - \M E^\top \widehat{\V x}_{k+1|n} \big).
\end{equation}

\bigskip 
We show below that the forward recursion is exactly what Kalman filtering does: at each step, it computes the least square estimate of the current state given the past information and the current measurement. 
We also show that the backward recursion allows to refine the Kalman estimates of the past states, an operation called \emph{smoothing}.

\subsection*{Forward recursion (Kalman filtering)}

Kalman filtering can be viewed as an implementation of the forward recursion used to solve
\eqref{eq:block_tridiag_improved}. In fact, substituting \eqref{eq:Eb1bkD} in \eqref{eq:Dkbk} and \eqref{eq:xnnGeneral},  and using Tool 2 and Tool 3, the forward recursion can be written as follows. 

Initialize
\begin{align*}
    \bar{\M D}_1& %= \M{\widetilde Q}_1^{-1}+\M H^\top \M R^{-1}\M H +\M F^\top \M Q^{-1}\M F 
    = \M P_1^{-1}+ \M F^\top \M Q^{-1}\M F, \\
    \bar{\V b}_1&= \M{\widetilde Q}_1^{-1}\M F\widehat{\V x}_0+\M H^\top\M R^{-1}\V z_1 = \M P_1^{-1} \bar{\V x}_1,
\end{align*}

for $k=2, 3, \dots,n-1$ compute % the recursion becomes
\begin{equation*}
\label{eq:DkbkKalman}
\begin{alignedat}{2}
%\M G_k &= \M E\,\bar{\M D}_{k-1}^{-1},\\
\bar{\M D}_k &%= \M D-\M G_k\M E^\top
%= \M D_k-\M E\,\bar{\M D}_{k-1}^{-1}\M E^\top 
= \M P_k^{-1} + \M F^\top \M Q^{-1}\M F,\\
\bar{\V b}_k &%= \V b_k-\M G_k\bar{\V b}_{k-1}
%= \V b_k-\M E\,\bar{\M D}_{k-1}^{-1}\bar{\V b}_{k-1} 
= \M{\widetilde Q}_k^{-1}\M F\bar{\V x}_{k-1}+\M H^\top\M R^{-1}\V z_k = \M P_k^{-1} \bar{\V x}_k 
\end{alignedat}
\end{equation*}

and %\eqref{eq:xnnGeneral} becomes
\begin{align}
    \bar{\M D}_n & %= \M D-\M G_k\M E^\top
%= \M D_n-\M E\,\bar{\M D}_{n-1}^{-1}\M E^\top 
%= (\M Q + \M F \M P_{n-1} \M F^\top)^{-1}+\M H^\top \M R^{-1}\M H 
= \M P_n^{-1} \qquad
\bar{\V b}_n %= \V b_k-\M G_k\bar{\V b}_{k-1}
%= \V b_n-\M E\,\bar{\M D}_{n-1}^{-1}\bar{\V b}_{n-1} 
= \M{\widetilde Q}_n^{-1}\M F\bar{\V x}_{n-1}+\M H^\top\M R^{-1}\V z_n \nonumber\\ %= \M P_n^{-1} \bar{\V x}_n.
 \widehat{\V x}_{n|n}&=\bar{\M D}_n^{-1}\,\bar{\V b}_n = \M P_n \big( \M{\widetilde Q}_n^{-1}\M F\bar{\V x}_{n-1}+\M H^\top\M R^{-1}\V z_n \big) \label{eq:xnnKalman}
\end{align}
where we have defined \,$\M{\widetilde Q}_k=\M Q + \M F \M P_{k-1} \M F^\top $,  $\M P_k = \big( \M{\widetilde Q}_k^{-1} + \M H^\top \M R^{-1}\M H \big)^{-1}$, %$= ((\M Q + \M F \M P_{k-1} \M F^\top)^{-1}+\M H^\top \M R^{-1}\M H)^{-1}$, 
and $\bar{\V x}_k = \M P_k \big( \M{\widetilde Q}_k^{-1} \M F \bar{\V x}_{k-1} + \M H^\top\M R^{-1}\V z_k \big)$, with $\bar{\V x}_0= \widehat{\V x}_0$. %We can see, comparing with  \eqref{eq:solution2}, that  $\bar{\V x}_k$ and $\M P_k$ are the output of the Kalman filter at time $k$. %In other words, 
We see that $\bar{\V x}_1$ and $\M P_1$ are the output \eqref{eq:solution2} of the Kalman filter given $\V z_1$, $\bar{\V x}_0$ and $\M P_0$; then $\bar{\V x}_2$ and $\M P_2$ are the output \eqref{eq:solution2} of the Kalman filter given $\V z_2$, $\bar{\V x}_1$ and $\M P_1$, and so on.

%At this point we have transformed the system into upper triangular form, so that 
%The least-squares solution \eqref{eq:xnnGeneral} for the last state is then %$\widehat{\V x}_n$
%\[
%\widehat{\V x}_{n|n}=\bar{\M D}_n^{-1}\,\bar{\V b}_n = \M P_n (\M{\widetilde Q}_n^{-1}\M F\bar{\V x}_{n-1}+\M H^\top\M R^{-1}\V z_n)\,.
%\]
Comparing with \eqref{eq:solution2}, we see that \eqref{eq:xnnKalman} coincides with the Kalman update at time $n$. 
%which, as anticipated, tells that the best estimate at time $n$ given $\V z_{1:n}$ is obtained by applying the Kalman recursion starting from the previous Kalman output. 
This proves that the Kalman filter equations yield, at time $n$, the least-squares/BLUE estimate of $\V x_n$
given the measurements $\V z_{1:n}$, by using as input the estimate from the previous step and proceeding recursively from the initial pair $\big( \widehat{\V x}_0,\M P_0 \big)$. Thus, the output of the Kalman filter at time $k$,
which we denoted by $\bar{\V x}_k$, is precisely the filtered estimate $\bar{\V x}_k=\widehat{\V x}_{k|k}$.

\bigskip

If we are interested in Problem~2, the past filtered estimates $\bar{\V x}_k=\widehat{\V x}_{k|k}$ can be refined by the backward substitution \eqref{eq:backwardgeneral}, which here reads as follows.

\noindent \emph{Backward recursion (smoothing).} For $k=n-1,n-2,\ldots,1$ compute
\begin{align}\label{eq:smoothing_elimination}
\widehat{\V x}_{k|n}
&= \big( \M P_{k}^{-1} + \M F^\top \M Q^{-1} \M F \big)^{-1} \big( \M P_{k}^{-1}\,\widehat{\V x}_{k|k} + \M F^\top \M Q^{-1}\,\widehat{\V x}_{k+1|n} \big) \nonumber \\
&= \M G_k\,\widehat{\V x}_{k+1|n} + \big( \M I-\M G_k \M F \big)\,\widehat{\V x}_{k|k} = \widehat{\V x}_{k|k} + \M G_k\big(\widehat{\V x}_{k+1|n}-\M F\widehat{\V x}_{k|k}\big),
\end{align}
where we used Tools 2) and 3), with $\M G_k \;\triangleq\; \M P_k \M F^\top \big( \M Q + \M F \M P_k \M F^\top \big)^{-1}$.

Equation \eqref{eq:smoothing_elimination} implements the Rauch--Tung--Striebel smoother for the state estimate \cite{Jaz:70,KaiSayHas:B00,Sim:06,GreAnd:15}, obtained here as the backward substitution associated with the block Gaussian elimination of \eqref{eq:block_tridiag_improved}.
For additional implementation aspects (e.g.\ covariance smoothing, numerical stability), the reader is referred to the referenced books \cite{Jaz:70,Kai:74, Har:90, Kay:93,KaiSayHas:B00, Sim:06, PetPetCam:09,GreAnd:15}.

\newpage
%-------------------------------------------------------
\noindent\emph{Smoothed error covariance.}
%-------------------------------------------------------
After running the Kalman filter forward we know, for each $k$, the filtered quantities $\widehat{\V x}_{k|k}$ and $\M P_k\equiv \M P_{k|k}$, obtained given measurements up to time $k$. However, the estimates can be improved using all the measurements--not only for the states as in \eqref{eq:smoothing_elimination}, but also for the covariances. In fact, looking at  \eqref{eq:block_tridiag_improved} we recognize that $\M M^{-1} = \big( \M A^\top \M R_{\V n}^{-1} \M A \big)^{-1}$ is the least square estimation error covariance of the stacked state vector. 
Therefore, the \emph{smoothed error covariances} $\M P_{k|n}$ are precisely the diagonal block entries of $\M M^{-1}$ 
\[
\M P_{k|n}\;\triangleq\;\mathbb E\Big\{ \big( \V x_k-\widehat{\V x}_{k|n} \big) \big( \V x_k-\widehat{\V x}_{k|n} \big)^\top \Big\} = \big( \M M^{-1} \big)_{kk}\,. 
\]
 
Similarly to the state smoothing recursion, these blocks can be computed recursively, running backward from $k=n-1$ down to $k=1$, as
\begin{equation}\label{eq:smoothing_cov}
\M P_{k|n} = \M P_k
+\M G_k \big( \M P_{k+1|n} - \widetilde{\M Q}_{k+1} \big) \,\M G_k^\top
\end{equation}
initialized with $\M P_{n|n}=\M P_n$. 

\medskip
The correction term
$\M G_k \big( \M P_{k+1|n} - \widetilde{\M Q}_{k+1} \big) \,\M G_k^\top$
propagates backward the extra information provided by future measurements, then typically $\M P_{k|n}\preceq \M P_k$, i.e., smoothing reduces the uncertainty also at time $k$.

%-------------------------------------------------------
Equation \eqref{eq:smoothing_cov} follows by the block Gaussian elimination step.  
%-------------------------------------------------------
In fact, the forward elimination recursion gives for $\M M$ the block bidiagonal factorization $\M M=\M L\,\M U$ with
\begin{equation}\label{eq:LU_blocks}
\M L=
\begin{bmatrix}
\M I & \M 0 & \cdots & \M 0\\
\M L_2 & \M I & \ddots & \vdots\\
\vdots & \ddots & \ddots & \M 0\\
\M 0 & \cdots & \M L_n & \M I
\end{bmatrix},
\qquad
\M U=
\begin{bmatrix}
\bar{\M D}_1 & \M E^\top & \cdots & \M 0\\
\M 0 & \bar{\M D}_2 & \ddots & \vdots\\
\vdots & \ddots & \ddots & \M E^\top\\
\M 0 & \cdots & \M 0 & \bar{\M D}_n
\end{bmatrix},
\qquad
\M L_k=\M E\,\bar{\M D}_{k-1}^{-1}.
\end{equation}

If we let $\M S=\M M^{-1}$ and denote by $\M S_{ii}$ its diagonal blocks, the 
backward recursion gives
\begin{align}\label{eq:Snn}
\M S_{nn}&=\bar{\M D}_n^{-1} 
\\ %\end{equation}
%then, for $k=n-1,n-2,\dots,1$,
%\begin{equation}\label{eq:diag_inv_recursion}
\M S_{kk}
&=
\bar{\M D}_k^{-1}
+
\bar{\M D}_k^{-1}\,\M E^\top\,\M S_{k+1,k+1}\,\M E\,\bar{\M D}_k^{-1} , \qquad k=n-1,n-2,\dots,1.
\end{align}
which, with \eqref{eq:Eb1bkD} and \eqref{eq:Dkbk}, gives \eqref{eq:smoothing_cov}.

\newpage

%-------------------------------------------------------
\section*{Worked numerical example (continued): smoothing}
%-------------------------------------------------------

We now exemplify the smoothing operation for the previous numerical example. 
There, we computed the output of the Kalman filter, i.e., the best estimate at each time $k$
using only the measurements up to time $k$:
$\widehat{\V x}_k=\widehat{\V x}_{k|k}$. 
Smoothing improves the estimates of past states (here $\V x_0$ and $\V x_1$) by exploiting
\emph{all} available measurements $z_1,z_2$.

\noindent
We use the smoother recursion \eqref{eq:smoothing_elimination} with smoothing gain $\M G_k \triangleq \M P_k \M F^\top \,\widetilde{\M Q}_{k+1}^{-1}$.

\paragraph{Smoothed state at $k=1$.}
Using $\M P_1, \M F^\top$ and $\widetilde{\M Q}_{2}$  from the worked example, we have
\[
\M G_1=
\begin{bmatrix}
0.754 & -0.368\\
0.100 & \phantom{-}0.936
\end{bmatrix}
\]
giving 
\[
\widehat{\V x}_{1|2}
=
\widehat{\V x}_1 + \M G_1 \big( \widehat{\V x}_{2|2} - \M F \widehat{\V x}_1 \big)
=
\begin{bmatrix}1.806\\3.923\end{bmatrix}
+
\begin{bmatrix}
0.754 & -0.368\\
0.100 & \phantom{-}0.936
\end{bmatrix}
\left(
\begin{bmatrix}4.231\\4.655\end{bmatrix}
-
\begin{bmatrix}3.768\\3.923\end{bmatrix}
\right)
\approx
\begin{bmatrix}1.886\\4.655\end{bmatrix}.
\]

\paragraph{Smoothed covariance at $k=1$.}
Using \eqref{eq:smoothing_cov} with $k=1$ and $\M P_{2|2}=\M P_2$,
\[
\M P_{1|2}
=
\M P_1 + \M G_1 \big( \M P_2 - \widetilde{\M Q}_2 \big) \M G_1^\top
\approx
\begin{bmatrix}
0.0318 & -0.0485\\
-0.0485 & 0.2345
\end{bmatrix}.
\]

\paragraph{Smoothed state at $k=0$.}
Using $\M P_0, \M F^\top$ and $\widetilde{\M Q}_{1}$ we get 
\[
\M G_0=
\begin{bmatrix}
0.988 & -0.489\\
0.005 & \phantom{-}0.988
\end{bmatrix}
\]
giving 
\[
\widehat{\V x}_{0|2}
=
\widehat{\V x}_0 + \M G_0 \big( \widehat{\V x}_{1|2} - \M F \widehat{\V x}_0 \big)
=
\begin{bmatrix}0\\4\end{bmatrix}
+
\begin{bmatrix}
0.988 & -0.489\\
0.005 & \phantom{-}0.988
\end{bmatrix}
\left(
\begin{bmatrix}1.886\\4.655\end{bmatrix}
-
\begin{bmatrix}2\\4\end{bmatrix}
\right)
\approx
\begin{bmatrix}-0.433\\4.646\end{bmatrix}.
\]

\paragraph{Smoothed covariance at $k=0$.}
Applying \eqref{eq:smoothing_cov} with $k=0$,
\[
\M P_{0|2}
=
\M P_0 + \M G_0 \big( \M P_{1|2} - \widetilde{\M Q}_1 \big) \M G_0^\top
\approx
\begin{bmatrix}
0.1462 & -0.1651\\
-0.1651 & 0.2382
\end{bmatrix}.
\]

\medskip
\noindent\emph{Comment.}
The smoother uses the \emph{future} measurement $z_2$ to refine earlier states: compared to
$\widehat{\V x}_{1|1}$ and $\widehat{\V x}_{0|0}$, the estimates $\widehat{\V x}_{1|2}$ and $\widehat{\V x}_{0|2}$
incorporate the information carried by the innovation at time $k=2$ propagated backward through the dynamics.
In this example, the update is especially visible on the velocity component, which is not measured directly but is inferred through the state model.

}

\section*{What we have learned?}

In this lecture note, we have derived the Kalman filter equations directly from the least squares estimation framework. We have shown that:

\begin{itemize}
\item  The Kalman filter is fundamentally the least squares estimator of the system state, with no need to assume Gaussian disturbances.
\item  The derivation follows naturally from the block-structured least squares problem combining state prediction and measurements.
\item  The same formulation accommodates correlations between process and measurement errors without requiring any structural changes to the equations.
\item  The optimality of the Kalman filter in the Gaussian case emerges as a special instance of the least squares interpretation.
\item  The sequential formulation, when the measurement noise covariance is diagonal, can be implemented using scalar updates, avoiding the explicit computation of matrix inverses and improving numerical efficiency.
{\color{revisionblue}
\item For estimating a state trajectory from noisy measurements, the least squares viewpoint leads directly to the Kalman filter and to backward recursions for smoothing.
}
\end{itemize}

Together, these points provide a unified and accessible view of the Kalman filter: it is not just an algorithm for combining predictions and measurements, but the natural solution to a least squares problem. This perspective simplifies the derivation, and clarifies why the standard form of the filter is both mathematically elegant and practically powerful.

\section*{Acknowledgment}
The fundamental research described in this paper was supported, in part, 
by the European Union through the Italian National Recovery and Resilience Plan of Next Generation EU  under Grant PE00000001-RESTART, 
by the National Research Foundation of Korea under Grant RS-2024-00409492, and 
by the Robert R. Taylor Professorship.
\newpage

%
%
% ---------------------------------------------------
% References
% ---------------------------------------------------
%
\bibliographystyle{IEEEtran}
\bibliography{Files/IEEEabrv,Files/StringDefinitions,Files/StringDefinitions2,Files/refs}
\vfill

\end{document}